\documentclass[12pt]{article}
\usepackage[latin1]{inputenc}
\usepackage{amsmath}
\usepackage{amsfonts}
\usepackage{amssymb}
\usepackage{color}
\usepackage{graphicx}
\newtheorem{theorem}{Theorem}[section]

\newtheorem{lemma}[theorem]{Lemma}

\newtheorem{remark}[theorem]{Remark}

\title{\bf  Backbone decomposition for continuous-state branching processes with immigration}
\author{{\sc A.E. Kyprianou\footnote{Department of Mathematical Sciences, University of Bath,
Claverton Down, Bath, BA2 7AY, U.K.}  \ \ and
  Y.-X.
Ren\footnote{LMAM School of Mathematical Sciences $\&$ Center for Statistical Science,
Peking University, Beijing 100871,
 P. R. China.} 
}
}

\begin{document}

\maketitle
\begin{abstract}

In the spirit of Duqesne and Winkel \cite{DW} and Berestycki et al. \cite{BKMS} we show that supercritical continuous-state branching process with a general branching mechanism and general immigration mechanism is equal in law to a continuous-time Galton Watson process with immigration with Poissonian dressing. The result also characterises the limiting backbone decomposition which is predictable from the work on consistent growth of Galton-Watson trees with  immigration in Cao and Winkel \cite{CW}.

\bigskip


\bigskip


\end{abstract}

\section{Introduction}\label{classicalbb}
In this article we are interested in the case that the $[0,\infty)$-valued strong Markov process with absorbing state at zero,  $X=\{X_t: t\geq 0\}$, is a conservative, supercritical continuous-state branching process with general branching mechanism
$\psi$ taking the form
\begin{equation*}
\psi(\lambda) = \alpha \lambda + \beta\lambda^2 + \int_{(0,\infty )} (e^{-\lambda x} - 1 + \lambda x \mathbf{1}_{\{x<1\}} )\Pi({\rm d}x),\,\, \lambda \geq 0,
\end{equation*}
where $\alpha\in\mathbb{R}$, $\beta\geq 0$ and $\Pi$ is a measure concentrated on $(0,\infty)$ which satisfies $\int_{(0,\infty)}(1\wedge x^2)\Pi({\rm d}x)<\infty$ and a general immigration mechanism $\varphi$ taking the form
\[
\varphi(\lambda) = \delta\lambda + \int_{(0,\infty)}(1- e^{-\lambda x})\nu({\rm d}x),
\]
where $\delta\geq 0$ and $\nu$ is a measure concentrated on $(0,\infty)$ which satisfies $\int_{(0,\infty)}(1\wedge x)\nu({\rm d}x)<\infty$. Our requirement that $X$ is supercritical and conservative means that we necessarily have that $\psi'(0+)<0$ and
\[
\int_{0+}\frac{1}{|\psi(\xi)|}{\rm d}\xi =\infty
\]
respectively.

The process $X$, henceforth denoted a $(\psi, \varphi)$-CSBP, can be described through its semi-group as follows. Suppose that  $\mathbb{P}_x$ denotes the law of $X$ on cadlag path space $D[0,\infty)$ when the process is issued from $x\geq 0$. Then the semi-group associated with the $(\psi,\varphi)$-CSBP can be described as follows. For all $x,\lambda \geq 0$ it necessarily follows that
\[
\mathbb{E}_x(e^{-\lambda X_t}) = e^{-  xu_t(\lambda) - \int_0^t \varphi(u_{t-s}(\lambda)) {\rm d}s}, \, \, t\geq 0,
\]
where $u_t(\lambda)$ uniquely solves the evolution equation
\begin{equation}
u_t(\lambda) + \int_0^t \psi(u_s(\lambda) ){\rm d}s = \lambda\label{int-equa}
\end{equation}
with initial condition $u_0(\lambda) = \lambda$.
Note in particular that $u_t(\lambda)$ describes the semi-group of the $(\psi, 0)$-CSBP.

Another process related to the $(\psi,0)$-CSBP is that of the $(\psi,0)$-CSBP conditioned to become extinguished. To understand what this means, let us momentarily recall that for all supercritical continuous-state branching processes (without immigration) the event $\{\lim_{t\uparrow\infty} X_t =0\}$ occurs with positive probability. Moreover, for all $x\geq 0$,
\[
\mathbb{P}_x(\lim_{t\uparrow\infty} X_t =0) = e^{-\lambda^* x}
\]
where $\lambda^*$ is the unique root on $(0,\infty)$ of the equation $\psi(\lambda) = 0$. Note that $\psi$ is strictly convex with the property that $\psi(0) = 0$ and $\psi(+\infty) = \infty$, thereby ensuring that the root $\lambda^*>0$ exists; see Chapter 8 and 9 of \cite{K} for further details. It is straightforward to show that
the law of $(X, \mathbb{P}_x)$ conditional on the event $\{\lim_{t\uparrow\infty} X_t =0\}$, say $\mathbb{P}^*_x$, agrees with the law of a $(\psi^*, 0)$-CSBP where
\begin{equation}\label{psi*}
\psi^*(\lambda ) = \psi(\lambda + \lambda^*).
\end{equation}
See for example Sheu \cite{Sheu1997}.

In Dusquene and Winkel \cite{DW} and Berestycki et al. \cite{BKMS} it was shown for the case that $\varphi\equiv 0$ that the law of process $X$ can be recovered from a supercritical continuous-time Galton-Watson process (GW), issued with a Poisson number of initial ancestors, and dressed in a Poissonian way using the law of the the original process conditioned to become extinguished.

To be more precise, they showed that for each $x\geq 0$, $(X, \mathbb{P}_x)$  has the same law as the  process $\{\Lambda_t : t \geq 0\}$ which has the following pathwise construction.
First sample from a continuous-time Galton-Watson process with branching generator
\begin{equation}
F(r) = q\left(\sum_{n\geq 0} p_n r^n - r\right) =  \frac{1}{\lambda^*}\psi(\lambda^*(1-r)).
\label{F}
\end{equation}
Note that in the above generator, we have that $q = \psi'(\lambda^*)$ is the rate at which individuals reproduce and $\{p_n: n\geq 0\}$ is the offspring distribution. With the particular branching generator given by \eqref{F}, $p_0 = p_1 =0$, and for $n\geq 2$,  $p_n : = p_n[0,\infty)$ where for $y\geq 0$,
\[
p_n({\rm d}y) =  \frac{1}{\lambda^* \psi'(\lambda^*)}\left\{\beta (\lambda^*)^2\delta_0({\rm d}y)\mathbf{1}_{\{n=2\}} + (\lambda^*)^n \frac{y^n}{n!} e^{-\lambda^*y} \Pi({\rm d}y)\right\}.
\]
If we denote the aforesaid GW process $Z = \{Z_t: t\geq 0\}$ then we shall also insist that $Z_0$ has a Poisson distribution with parameter $\lambda^*x$.
Next, {\it dress} the life-lengths of $Z$ in such a way that a $(\psi^*,0)$-CSBP is independently grafted on to each edge of $Z$ at time $t$ with rate
\[
2\beta {\rm d}\mathbb{N}^* + \int_0^\infty y e^{-\lambda^* y}\Pi({\rm d}y){\rm d}\mathbb{P}^*_y.
\]
Here the measure $\mathbb{N}^*$ is the Dykin-Kuznetsov excursion measure on the space $D[0,\infty)$ which satisfies
\[
\mathbb{N}^*(1- e^{-\lambda X_t}) = u^*_t(\lambda)=- \frac{1}{x}\log\mathbb{E}^*_x(e^{-\lambda X_t})
\]
for $\lambda,t\geq 0$, where $u^*_t(\lambda)$ is the unique solution to the integral equation
\begin{equation}\label{int-equa2}
u^*_t(\lambda) + \int_0^t \psi^*(u^*_s(\lambda)) = \lambda,
\end{equation}
with initial condition $u^*_0(\lambda) = \lambda$.
See \cite{Dynkin-Kuznetsov} for further details. Moreover, on the event that an individual dies and branches into $n\geq 2$ offspring, with probability $p_n({\rm d}x)$, an additional independent $(\psi^*, 0)$-CSBP is grafted on to the branching point with initial mass $x\geq 0$. The quantity $\Lambda_t$ is now understood to be the total dressed mass present at time $t$ together with the mass present at time $t$ in an independent $(\psi^*,0)$-CSBP issued at time zero with initial mass $x$.

\bigskip

 Our objective here is to describe a similar decomposition for the $(\psi, \varphi)$-CSBP. In the case that we include immigration, it will turn out that the backbone is rather naturally replaced by a continuous-time Galton-Watson process with immigration.

\section{Backbone decomposition}

In order to describe the backbone decomposition for the $(\psi, \varphi)$-CSBP, let us first remind ourselves of the basic structure of a continuous-time Galton-Watson process with immigration. Such processes are characterised by the two generators $(F, G)$ where, as mentioned before,
\[
F(r) = q\left(\sum_{n\geq 0} p_n r^n - r\right)
\]
encodes the fact that individuals live for an independent and exponentially distributed length of time, after which they give birth to a random number of offspring with distribution $\{p_n: n\geq 0\}$, and
\[
G(r) = p\sum_{n\geq 0} \pi_n r^n,
\]
 reflecting the fact that at times of a Poisson arrival process with rate $p>0$, a random number of immigrants with distribution $\{\pi_n: n\geq 0\}$ issue independent copies of a continuous-time Galton-Watson process with generator $F$.

Our forthcoming backbone decomposition will be built from an $(F,G)$-GW process with
$F$ given by (\ref{F}) and
\begin{equation}
G(r) = \varphi(\lambda^*) - \varphi(\lambda^*(1-r))
\label{G}
\end{equation}
It can be seen from the above expression for $G(r)$ that $p= \varphi(\lambda^*)$. To describe the distribution  $\{\pi_n: n\geq 0\}$ let us introduce an associated probability measure,
concentrated on $\{1,2,\cdots\}\times(0,\infty)$,
{\color{black}
\begin{equation}
\pi_n({\rm d}y) = \frac{1}{\varphi(\lambda^*)} \left[ (\delta\lambda^*)\delta_0({\rm d}y)\mathbf{1}_{(n=1)}+ \frac{(\lambda^*y)^n}{n!} e^{-\lambda^*y} \nu({\rm d}y)\right].
\label{imm}
\end{equation}
}
It is straightforward to check  that, in (\ref{G}), $\pi_0:=0$,
 $\pi_n : = \pi_n(0,\infty), \, n\geq 1$ and $p=\varphi(\lambda^*)$ respectively.

\bigskip

Fix $x>0$. Our backbone decomposition for the process $(X, \mathbb{P}_x)$ will consist of the bivariate Markov process $(Z, \Lambda) = \{(Z_t , \Lambda_t): t\geq 0\}$ valued in $\{0,1,2,\ldots\}\times [0,\infty)$. Here the process, $Z$, the backbone, is an $(F,G)$-GW process as described above with the additional property that $Z_0$ is Poisson distributed in number with rate $\lambda^* x$. The process of continuous mass, $\Lambda$, is described as follows.

\bigskip

(i) As in \cite{BKMS}, along the life length of each individual alive in the process $Z$, there is Poissonian dressing with rate
\begin{equation}
2\beta {\rm d}\mathbb{N}^* + \int_0^\infty y e^{-\lambda^* y}\Pi({\rm d}y){\rm d}\mathbb{P}^*_y.
\label{dressing}
\end{equation}

\medskip

(ii) At the branch points of $Z$, on the event that there are $n$ offspring, an additional copy of a $(\psi^*, 0)$-CSBP
with initial mass $y\geq 0$ is issued with probability $p_n({\rm d}y)$.

\medskip

(iii) At the same time, along the time-line between  each immigration of $Z$, there is again Poissonian dressing with rate
\begin{equation}
\delta{\rm d}\mathbb{N}^* + \int_0^\infty  e^{-\lambda^* y}\nu({\rm d}y){\rm d}\mathbb{P}_y^*.
\end{equation}

\medskip

(iv) Moreover, on the event that there are $n\geq 1$ immigrants in $Z$, an additional copy of a $(\psi^*,0)$-CSBP with initial mass $y\geq 0$ is issued with probability $\pi_n({\rm d}y)$.

\bigskip

\begin{figure}[t]
  \begin{center}
    \includegraphics[width=13cm]{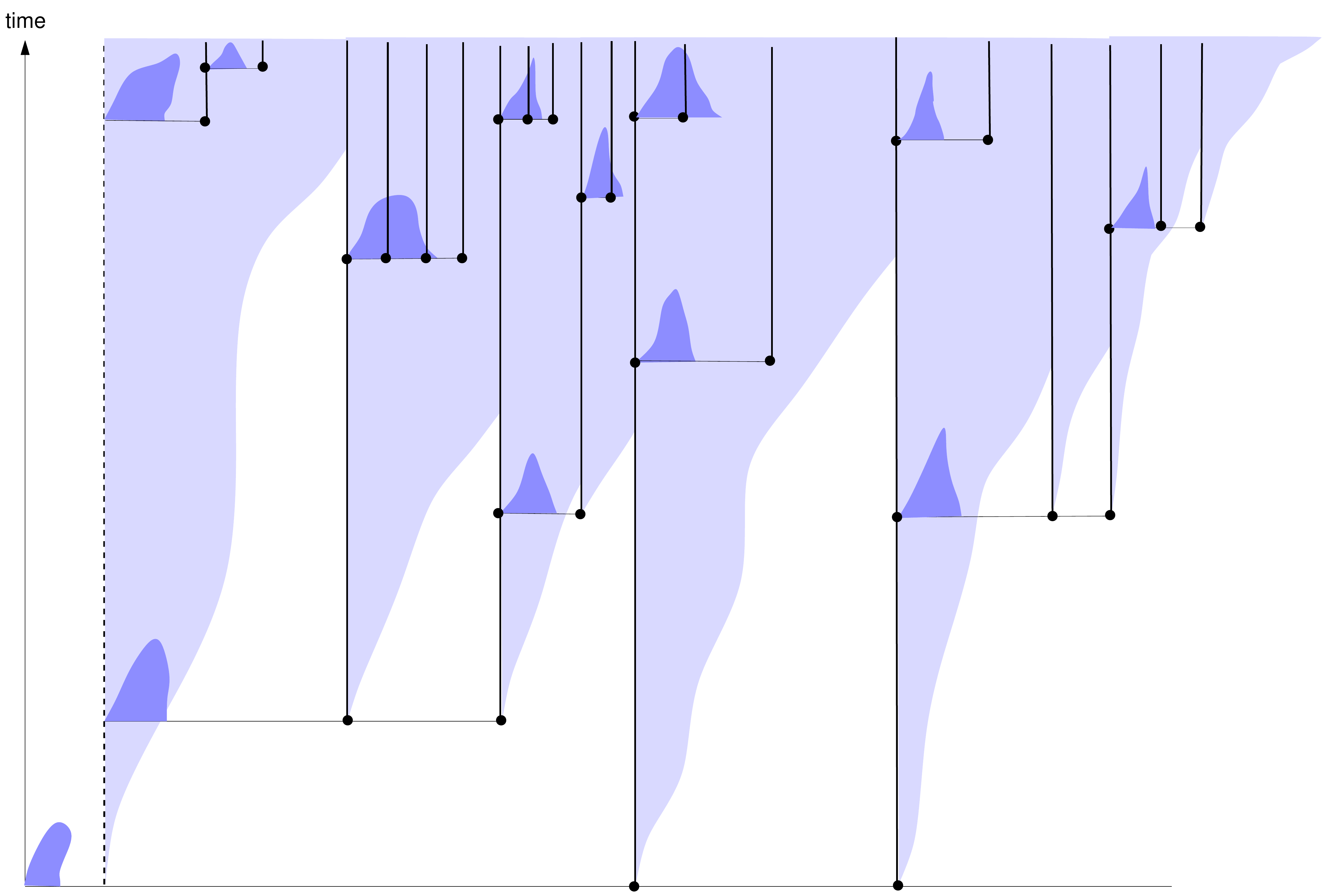}
  \end{center}
  \caption{The diagram above gives a symbolic representation of the backbone decomposition for the $(\psi, \varphi)$-CSBP. Working from left to right: An independent copy of a $(\psi^*,0)$-CSBP (shaded dark) is issued at time zero with initial mass $x$ together with an $(F,G)$-GW process which admits a Poisson distributed number of initial individuals with rate $\lambda^*x$. Along the (vertical dotted) time-line of the immigration process the dressing  (shaded light) has rate $\delta{\rm d}\mathbb{N}^* + \int_0^\infty  e^{-\lambda^* y}\nu({\rm d}y){\rm d}\mathbb{P}_y^*$ and additional independent $(\psi^*,0)$-CSBPs  (shaded dark) are grafted on  at times of immigration such that the probability there are $n$ simultaneous immigrants with grafted mass of initial size $y\geq 0$ is  $\pi_n({\rm d}y)$. Along the life length of individuals in the $(G,F)$-GW process (vertical black lines)  there is dressing (shaded light) at rate ${\color{black}2}\beta {\rm d}\mathbb{N}^* +
  \int_0^\infty y e^{-\lambda^* y}\Pi({\rm d}y){\rm d}\mathbb{P}^*_y$ with additional independent mass (shaded dark) grafted on at branching times such that the probability of there {\color{black}being $n$ offspring with grafted mass of initial size $y>0$} is $p_n({\rm d}y)$.}
  \label{F1}
\end{figure}

The quantity $\Lambda_t$ is now taken to be the total dressed mass present at time $t$ together with the mass at time $t$ in an independent $(\psi^*,0)$-CSBP issued at time zero with initial mass $x$. Figure \ref{F1} gives a pictorial representation of this decomposition. Henceforth we shall denote the law of the process $(Z,\Lambda)$ by $\mathbf{P}_x$.

 \begin{theorem}[Backbone decomposition for $(\psi, \varphi)$-CSBP]\label{bb}
 Fix $x>0$.
The law of $(X, \mathbb{P}_x)$ agrees with that of $(\Lambda, \mathbf{P}_x)$.
  Moreover, for all $t\geq 0$, the law of $Z_t$ given $\Lambda_t$ is that of a Poisson random variable with law $\lambda^* \Lambda_t$.
 \end{theorem}

\begin{remark}\rm  The above decomposition complements the recent work of Cao and Winkel \cite{CW}. In their paper, it is shown how to consistently grow GW trees with immigration in such a way that, with suitable rescaling, the resulting total mass at each fixed time converges in law to that of a $(\psi, \varphi)$-CSBP process. In some sense, the decomposition in Theorem \ref{bb} gives a slightly richer description of what the rescaled GW trees with immigration in \cite{CW} will converge to.
\end{remark}

\begin{remark}\rm Before progressing to the proof, we note that the above theorem can also be cited in the setting of a general superprocess where the motion, taken as a  general Borel right Markov process with Lusin state space, is independent of the branching mechanism (now reading $Z, X$ and $\Lambda$ as random measures) with minor modification to the forthcoming proof, providing one insists further that $|\psi'(0+)|<\infty$. The additional condition  is inherited from Berestycki et al. \cite{BKMS}.  Whilst this condition is not required in the case that motion is neglected, \cite{BKMS} require it as soon as spatial considerations come into play.
\end{remark}

 \section{Proof of  main result}

 We first need a result in \cite{BKMS}  which was originally stated for superprocesses. We use it here in a reduced form (the spacial movement of particles is reduced to a fixed point).

 \begin{lemma}\label{lemma} Let $(Z^\emptyset, \Lambda^\emptyset)$ be a copy of the backbone decomposition for a $(\psi, 0)$-CSBP, where  the process $Z^\emptyset$, the backbone, is an $(F,0)$-GW process as described above with the additional property that  $Z^\emptyset_0=n\in\{0,1,2,\ldots\}$, the process of continuous mass, $\Lambda^\emptyset$, is described as above with  the additional property that $\Lambda^\emptyset_0=y$.  Let ${\bf P}^\emptyset_{(y,n)}$ be the law of $(Z^\emptyset, \Lambda^\emptyset)$. Then \[
 {\bf E}^\emptyset_{(y,n)}(r^{Z^\emptyset_t} e^{-\theta \Lambda^\emptyset_t}) = e^{-yu^*_{t}(\theta)-nw_{t}(r,\theta)},
 \] where \begin{equation}
\lambda^*(1 - e^{ - w_t(r, \theta)}) = u_t(\theta + \lambda^*(1-r)) - u_t^*(\theta).
\label{w}
 \end{equation}
 \end{lemma}

 {\bf Proof:}\quad According to Theorem 1 in \cite{BKMS},
 \[{\bf E}^\emptyset_{(y,n)}(r^{Z^\emptyset_t}e^{-\theta\Lambda^\emptyset_t})=e^{-yu^*_{t}(\theta)-nw_{t}(r,\theta)},\]
 where $e^{-w_{t}(r,\theta)}$ is the unique $[0,1]$-valued solution to the integral equation
 \[e^{-w_{t}(r,\theta)}=r+\frac{1}{\lambda^*}\int^t_0{\rm d}s[\psi^*(-\lambda^*e^{-w_{t-s}(r,\theta)}+u^*_{t-s}(\theta))-\psi^*(u^*_{t-s}(\theta))]
 \]
 for $t\ge 0.$ With the help of \eqref{psi*} and \eqref{int-equa2}, it is  straightforward to show that
 $u_t^*(\theta)+\lambda^*(1 - e^{ - w_t(r, \theta)})$ solves \eqref{int-equa} with initial condition $\theta+\lambda^*(1-r)$. Therefore we have \[
\lambda^*(1 - e^{ - w_t(r, \theta)}) = u_t(\theta + \lambda^*(1-r)) - u_t^*(\theta)
\]
as required.
\hfill$\square$

\bigskip

 {\bf Proof of Theorem  \ref{bb}:}\quad For the first part we need to show that the process $(\Lambda, \mathbf{P}_x)$ is Markovian and its semi-group agrees with that of $(X, \mathbb{P}_x)$. For the second part it suffices to show that for $r\in[0,1]$ and $\theta\geq 0$,
 \begin{equation}
 \mathbf{E}_x( r^{Z_t}e^{-\theta \Lambda_t }) = \mathbf{E}_x(e^{-(\theta + \lambda^*(1-r))\Lambda_t }).
 \label{*}
 \end{equation}
 It fact, a little thought shows that both of these facts can be simultaneously established by proving that for all $x\geq 0$, $r\in[0,1]$ and $\theta\geq 0$,
 \begin{equation}
 \mathbf{E}_x(r^{Z_t}e^{-\theta \Lambda_t }) = e^{- xu_t(\theta + \lambda^*(1-r))}.
 \label{**}
 \end{equation}
 Indeed, note that (\ref{**}) directly implies (\ref{*}) and by setting $r=1$ in (\ref{*}) we also see that $\Lambda$ has the required semi-group.


 To this end, let us split the process $(Z, \Lambda)$ in to the independent sum of processes $(Z^\emptyset, \Lambda^\emptyset)$ and $(Z^I, \Lambda^I)$ where the first is an independent copy of the backbone decomposition for a $(\psi, 0)$-CSBP and $(Z^I, \Lambda^I)$ is the part of $Z$  rooted at immigration times together with its dressing. Note immediately by independence  we have that
 \[
 \mathbf{E}_x(r^{Z_t}e^{-\theta \Lambda_t }) = \mathbf{E}_x(r^{Z^\emptyset_t}e^{-\theta \Lambda^\emptyset_t }) \mathbf{E}_x(r^{Z^I_t}e^{-\theta \Lambda^I_t }) =  e^{-x u_t (\theta + \lambda^*(1-r))} \mathbf{E}_x(r^{Z^I_t}e^{-\theta \Lambda^I_t }),
 \]
 where the second equality follows from the Poissonization that is known to hold for the backbone embedding of $(\psi,0)$-CSBPs as described in \cite{BKMS} (see also the discussion in Section \ref{classicalbb}).

It therefore suffices to prove that for all $x\geq 0$, $s\in[0,1]$ and $\theta\geq 0$
\[
 \mathbf{E}_x(r^{Z^I_t}e^{-\theta \Lambda^I_t }) = e^{- \int_0^t \varphi(u_{t-s} (\theta + \lambda^* (1-r))) {\rm d}s}.
\]

With this as our goal, let us now write for each $t\geq 0$,
\[
\Lambda^I_t = \Lambda^{I,1}_t + \Lambda^{I,2}_t,
\]
 where $\Lambda^{I,1}_t$ is the mass at time $t$ due to the Poissonian dressing along  the time-line between each immigration of $Z$  and $\Lambda^{I, 2}$ is the mass at time $t$ due to the dressing  at immigration times together with the dressing of the immigrating $(F, 0)$-GW processes. First note that with the help of Campbell's Formula,
 \begin{eqnarray}
 \mathbf{E}_x(e^{-\theta \Lambda^{I,1}_t}) &=& \exp\left\{
 - \int_0^t  {\rm d}s\cdot \delta\mathbb{N}^*(1 - e^{-\theta X_{t-s}})  {\color{black}-} \int_{(0,\infty)} e^{-\lambda^* y}\nu({\rm d}y)\mathbb{E}^*_y(1 - e^{-\theta X_{t-s}})\right\}\notag\\
 &=&\exp\left\{
 -\int_0^t{\rm d}s\cdot \delta u^*_{t-s}(\theta) {\color{black}-} \int_{(0,\infty)}(1- e^{-y u^*_{t-s} (\theta)}) e^{-\lambda^* y}\nu({\rm d}y)
 \right\}\notag\\
 &=& \exp\left\{ - \int_0^t {\rm d}s\cdot\varphi^*(u^*_{t-s}(\theta))  \right\},
 \label{1}
  \end{eqnarray}
  where $$\varphi^*(\lambda) := \varphi(\lambda +\lambda^* ) - \varphi(\lambda^*)=\delta \lambda+\int^\infty_0(1-e^{-\lambda y})e^{-\lambda^* y}\nu(dy).$$
Recalling that the  immigration of $Z$ is characterised by the $G$, by using Lemma \ref{lemma} and applying  Campbell's Formula, we have
 \begin{eqnarray}
\lefteqn{ \mathbf{E}_x(r^{Z^I_t} e^{-\theta {\color{black}\Lambda^{I,2}_t}})}&&\notag\\
 &=&
 \exp\left\{
 - \int_0^t {\rm d}s\cdot\varphi(\lambda^*)\sum_{n\geq 1}\int_0^\infty {\color{black}\pi_n({\rm d}y)}{\color{black}(1 - e^{-y u^*_{t-s}(\theta) -n w_{t-s}(r, \theta)})}
 \right\}\notag\\
 &=&
 \exp\Bigg\{
 -\int_0^t {\rm d}s\cdot\Bigg({\color{black}\delta \lambda^*}+\int_{(0,\infty)}(1- e^{-\lambda^* y})\nu({\rm d}y)\notag\\
 &&\hspace{2cm}- \int_{(0,\infty)}\sum_{n\geq 1}\frac{(\lambda^* y e^{-w_{t-s}(r,\theta)})^n}{n!} e^{-\lambda^* y} e^{-y u^*_{t-s}(\theta) }\nu({\rm d}y){\color{black}-\delta\lambda^*e^{-w_{t-s}(r, \theta)}}
 \Bigg)\Bigg\}\notag\\
 &=&
\exp\Bigg\{
 -\int_0^t{\rm d}s\cdot
\Bigg(\varphi(\lambda^*)
- \int_{(0,\infty)}(\exp\{ \lambda^* y e^{-w_{t-s}(r, \theta)}  \} -1)
e^{-y(\lambda^* + u^*_{t-s}(\theta))}\nu({\rm d}y)\notag\\
&&\hspace{11cm}-\delta\lambda^*e^{-w_{t-s}(r, \theta)}\Bigg)
 \Bigg\}\notag
 \\
  &=&
 \exp\left\{
 -\int_0^t{\rm d}s\cdot
 \left(\varphi(\lambda^*) + \varphi^*_{u^*_{t-s}(\theta)}(-\lambda^* e^{-w_{t-s} (r,\theta)})\right)
\right\},
\label{2}
 \end{eqnarray}
 where for $u\geq -\lambda^*$,
 \begin{eqnarray*}
 {\color{black}\varphi^*_{u}}(\lambda ) &= &\varphi^*(\lambda + u) - \varphi^*(u) =
 \varphi(\lambda + \lambda^* + u) - \varphi(\lambda^* + u)\\
 &=&\delta\lambda+\int_{(0,\infty)}(1-e^{-\lambda y})e^{-y(\lambda^*+u)}\nu({\rm d}y).
\end{eqnarray*}
 Putting the pieces together in (\ref{1}) and (\ref{2}) with the help of (\ref{w}), we see that
 \begin{eqnarray*}
 \mathbf{E}_x(r^{Z^I_t}e^{-\theta \Lambda^I_t })
& = &
\exp\left\{-
\int_0^t{\rm d}s\cdot \varphi(u^*_{t-s}(\theta) +\lambda^*(1- e^{-w_{t-s}(r,\theta)}))
\right\}\\
&=&\exp\left\{-
\int_0^t{\rm d}s\cdot
\varphi(u_{t-s}(\theta + \lambda^*(1-r))
\right\}
 \end{eqnarray*}
  as required.\hfill$\square$

\section*{Acknowledgements}
 The research of YXR  is supported in part by NNSF of China
 (Grant No. 10871103 and 10971003) and Specialized Research Fund for the Doctoral Program of Higher Education.


\begin{thebibliography}{99}

\bibitem{BKMS} J. Berestycki, A.E. Kyprianou and A. Murillo-Salas (2011): The prolific backbone for supercritical superprocesses. {\it Stoch. Proc. Appl.} {\bf 121}, 1315-1331.


\bibitem{CW} Cao, X. and Winkel, M. (2010): Growth of Galton-Watson trees: immigration and lifetimes. \texttt{arXiv:1004.3061v1 [math.PR]}


 \bibitem{DW} T. Duquesne and M. Winkel (2007): Growth of L\'evy trees. {\it Probab. Theory  Relat. Fields}, {\bf 139}, 313--371.



 \bibitem{Dynkin-Kuznetsov} E.B. Dynkin and S.E. Kuznetsov (2004):  $\mathbb{N}$-measures for branching Markov exit systems and their applications to differential equations. {\it Probab. Theory Relat. Fields.} {\bf 130}, 135-150.

\bibitem{K} A.E. Kyprianou (2006): {\it Introductory lectures on fluctuations of L\'evy processes with applications.} Springer.


\bibitem{Sheu1997}  Y.C. Sheu (1997):  Lifetime and compactness of range for $\psi$-super-Brownian motion with a general
                branching mechanism.
              {\em Stoch. Proc. Appl.} {\bf 70},
              129-141.


\end{thebibliography}
\end{document}